\documentclass[12pt,psamsfonts]{amsart}
\usepackage{latexsym}
\usepackage{amsmath, amsfonts, amssymb, amscd}
\usepackage[height=22.5cm, width=16cm]{geometry}
\usepackage[colorlinks=true, pdfstartview=FitV, linkcolor=blue,
citecolor=blue, urlcolor=blue]{hyperref}

\DeclareMathOperator{\grad}{grad} \DeclareMathOperator{\id}{id}
\DeclareMathOperator{\Ad}{Ad} \DeclareMathOperator{\ad}{ad}

\newtheorem{definition}{Definition}[section]
\newtheorem{theorem}[definition]{Theorem}
\newtheorem{lemma}[definition]{Lemma}
\newtheorem{proposition}[definition]{Proposition}
\newtheorem*{theorem11}{Theorem 1.1}

\theoremstyle{remark}
\newtheorem*{remark}{Remark}

\title{Semistable points with respect to real forms}
\author{Peter Heinzner}

\address{Fakult\"at f\"ur Mathematik\\
  Ruhr Universit\"at Bochum\\
  Universit\"atsstrasse 150\\
  D - 44780 Bochum}
\email{heinzner@cplx.rub.de}\email{henrik@cplx.rub.de}
\author{Henrik St\"otzel}
\thanks{The authors are partially supported by the
  Sonderforschungsbereich SFB/TR12 of the Deutsche
  Forschungsgemeinschaft}

\begin{document}

\begin{abstract}
We consider actions of real Lie subgroups $G$ of complex reductive
Lie groups on K\"ahlerian spaces. Our main result is the openness
of   the set of semistable points with respect to a momentum map
and the action of $G$.
\end{abstract}

\maketitle

\section{Introduction}

Throughout this article we consider a compact Lie group $U$ and a
Lie subgroup $G$ of the complexified group $U^\mathbb{C}$. The
corresponding Lie algebras will be denoted by German letters. We
will assume that $G$ is compatible with the Cartan decomposition
$U^\mathbb{C}=U\exp\,i\mathfrak{u}$ of $U^\mathbb{C}$. This means
that we have a diffeomorphism $K\times\mathfrak{p}\to G$,
$(k,\xi)\mapsto k\exp\,\xi$, where $K:=G\cap U$ and
$\mathfrak{p}:=\mathfrak{g}\cap i\mathfrak{u}$.

Let $Z$ be a complex space endowed with a holomorphic action of
$U^\mathbb{C}$ and a $U$-invariant K\"ahlerian structure $\omega$.
For simplicity we will assume that $\omega$ is smooth when
restricted to any smooth submanifold of $Z$. We also assume that
there is a continuous $U$-equivariant momentum map $\mu\colon Z\to
\mathfrak{u}^*$ which is smooth on every smooth submanifold of
$Z$. In particular, this means that the restriction of $\omega$ to
a smooth $U^\mathbb{C}$-stable submanifold $M$ of $Z$ is given by
a K\"ahlerian form in the usual sense. Hence $M$ is a Riemannian
manifold. For $\xi\in\mathfrak{u}$ we have the function
$\mu^\xi\colon M\to\mathbb{R}$, $\mu^\xi(z):=\mu(z)(\xi)$, with
gradient
$\grad(\mu^\xi)(z)=\left.\frac{d}{dt}\right|_0\exp(it\xi)\cdot
z=(i\xi)_Z(z)=J\xi_Z(z)$. Here $J$ denotes the complex structure
on the tangent bundle $TM$ and $\xi_Z$ is the vector field on $Z$
corresponding to the action of the one parameter group  $t\mapsto
\exp (t\xi)$.

For a subspace $\mathfrak{m}$ of $\mathfrak{u}$ we have the map
$\mu_\mathfrak{m}\colon Z\to \mathfrak{m}^*$ which is given by
composing $\mu$ with the adjoint of the inclusion map
$\mathfrak{m}\to\mathfrak{u}$. We call $\mu_{\mathfrak{m}}$ the
restriction of $\mu$ to $\mathfrak{m}$. For a subset $Q$ of $Z$ we
set $\mathcal{S}_G(Q):=\{z\in Z;\ \overline{G\cdot z}\cap
Q\ne\emptyset\}$. Here $\overline{G\cdot z}$ denotes the closure
of $G \cdot z$ in $Z$. For $\beta \in i\mathfrak{p}^*$ let
$\mathcal
M_{i\mathfrak{p}}(\beta):=\mu_{i\mathfrak{p}}^{-1}(\beta)$ and
since $\beta=0$ plays a prominent role we set $\mathcal
M_{i\mathfrak{p}}:=\mathcal M_{i\mathfrak{p}}(0)$.

The \emph{set of semi-stable points} of $Z$ with respect to
$\mu_{i\mathfrak{p}}$ and the $G$-action is by definition the set
$\mathcal{S}_G(\mathcal M_{i\mathfrak{p}})$. In \cite{HS05} it is
shown that $\mathcal{S}_G(\mathcal M_{i\mathfrak{p}})$ is open
when $U$ is abelian. Here we give a proof of a more general
result.

\begin{theorem} \label{main}
For every $\beta\in i\mathfrak{p}^*$ the set
$\mathcal{S}_G(\mathcal M_{i\mathfrak{p}}(\beta))$ is open in $Z$.
\end{theorem}

The authors would like to thank Christian Miebach and Patrick
Sch\"utzdeller for  several interesting discussions on the
material of this paper.

\section{Outline of the proof}

First we consider the case $\beta=0$. It is sufficient to show
that every point $x\in\mathcal M_{i\mathfrak{p}}$ has an open
neighborhood $V(x)$ of $x$ in $Z$ such that $V(x)\subset
\mathcal{S}_G(\mathcal M_{i\mathfrak{p}})$, since in this case
$\mathcal{S}_G(\mathcal M_{i\mathfrak{p}})=G\cdot V$ is open in
$Z$ where $V$ denotes the union of the sets $V(x)$.

Let $\lVert\ \rVert\colon\mathfrak{u}^*\to\mathbb{R}$ denote a
$U$-invariant norm function and let $\eta_{i\mathfrak{p}}\colon
Z\to\mathbb{R}$, $\eta_{i\mathfrak{p}}(z)=\frac{1}{2}\lVert
\mu_{i\mathfrak{p}}(z)\rVert^2$. Using convexity properties of the
momentum map with respect to special commutative subgroups (see
section~\ref{commutative}\,\!), we will show the existence of
$G$-stable neighborhoods $\Omega$ of points in $\mathcal
M_{i\mathfrak{p}}$ such that for some $r>0$ the sets
$V:=\Omega\cap\eta_{i\mathfrak{p}}^{-1}([0,r))$ are relatively
compact in $Z$.

For simplicity assume now that $Z$ is smooth. The flow $\psi$ of
the gradient vector field $\grad (\eta_{i\mathfrak{p}})$ is
tangent to the $G$-orbits. Moreover for $z\in V$ the integral
curve $t\mapsto \psi_z(t)$ exists for all negative $t$ and is
contained in $V$. If $t_n\to -\infty$, then possibly after going
over to a subsequence a limit point $z_0=\lim \psi_z(t_n)$ exists.
We have $z_0\in \overline{V\cap G\cdot z}$ and $z_0$ is a critical
point of $\eta_{i\mathfrak{p}}$. It then remains to show that $V$
can be chosen such that $\mathcal M_{i\mathfrak{p}}\cap
\overline{V}$ equals the set of critical point of
$\eta_{i\mathfrak{p}}$ which are contained in $\overline{V}$.

The case where $\beta$ is arbitrary is done by shifting the
momentum map. The usual procedure is to endow the coadjoint orbit
$O:=U\cdot\beta$ with the structure of a K\"ahlerian manifold such
that $\alpha\mapsto -\alpha$ defines a momentum map on $O$. Since
$O$ is compact the $U$-action extends to a holomorphic
$U^\mathbb{C}$-action $(g,\alpha)\mapsto g\cdot \alpha$, i.e.,
$O=U\cdot \beta=U^\mathbb{C}\cdot\beta$ holds. Using the shifting
procedure means to replace $Z$ by $Z\times O$ on which
$(z,\alpha)\mapsto \mu(z)-\alpha$ defines a momentum map $\hat\mu$
and study the set of semistable points with respect to $\hat\mu$
as well as its relation to $\mathcal{S}_G(\mathcal
M_{i\mathfrak{p}}(\beta))$. In fact in the case where
$G=U^\mathbb{C}$ this gives rather directly Theorem~\ref{main}. In
the general case one has to take into account the basic
observation that $G\cdot\beta=K\cdot\beta$ holds for $\beta\in
i\mathfrak{p}^*$.

\section{The structure of compatible subgroups}\label{compatible}

As before let $U$ be a compact Lie group, $U^\mathbb{C}$ its
complexification and let $G=K\exp\mathfrak{p}$ be a compatible
subgroup of $U^\mathbb{C}$.

If $G=K\exp(\mathfrak{p})$ is not closed in $U^\mathbb{C}$, then
its topological closure $\overline G=\overline
K\exp(\mathfrak{p})$ is a closed Lie-subgroup of $U^\mathbb{C}$
which is compatible with the Cartan decomposition of
$U^\mathbb{C}$. Moreover since $\overline{Gz}=\overline{\overline
Gz}$ we have $\mathcal{S}_G(\mathcal
M_{i\mathfrak{p}}(\beta))=\mathcal{S}_{\overline G}(\mathcal
M_{i\mathfrak{p}}(\beta))$ for $\beta\in i\mathfrak{p}^*$. This
implies that for the proof of Theorem~\ref{main} we may assume
that $G$ is a compatible closed subgroup of $U^\mathbb{C}$. The
main advantage of $G$ being closed is that $K$ is a compact
subgroup of $U$ and that for a maximal subalgebra $\mathfrak{a}$
of $\mathfrak{p}$ we have $G=KAK$ for $A:=\exp \mathfrak{a}$.

The remaining part of this section is not used until
section~\ref{shifting} where we need precise information about the
coadjoint action of $U$. The main input there is the fact that the
coadjoint orbit of $U$ through a point $\beta\in i\mathfrak{p}^*$
is in a natural way a complex $U^\mathbb{C}$-homogeneous manifold
$O$ such that the orbit $G\cdot \beta$ equals $K\cdot\beta$.

To begin with recall that the group $U^\mathbb{C}$ has a unique
structure of an affine algebraic group. Its algebra of regular
functions consists of the $U$-finite holomorphic functions on
$U^\mathbb{C}$. If $G$ is any subgroup of $U^\mathbb{C}$, then its
algebraic (respectively analytic) Zariski closure is the smallest
closed algebraic (respectively analytic) subvariety of
$U^\mathbb{C}$ which contains $G$. The Zariski closure is
automatically an algebraic (respectively analytic) subgroup of
$U^\mathbb{C}$. Moreover in the compatible case the notion of
algebraic and analytic Zariski closure of $G$ coincide. The
simplest way to see this is to use the fact that any closed
complex subgroup of $U^\mathbb{C}$ which is stable with respect to
the Cartan involution $\Theta\colon U^\mathbb{C}\to U^\mathbb{C}$,
$\Theta(u\exp i\xi):=u\exp(-i\xi)$ where $u\in U$ and
$\xi\in\mathfrak{u}$ is the complexification of its $\Theta$-fixed
points. A slightly more precise statement is the following simple

\begin{lemma}\label{Zariskilemma}
Let $G=K\exp\mathfrak{p}$ be a connected compatible Lie subgroup
of $U^\mathbb{C}$ and let $U_0$ be the smallest closed subgroup of
$U$ which contains $\exp(\mathfrak{k}+i\mathfrak{p})$. Then the
Zariski closure of $G$ in $U^\mathbb{C}$ equals
$(U_0)^\mathbb{C}=U_0\exp(i\mathfrak{u}_0)$. In particular $G$ is
compatible with the Cartan decomposition of $(U_0)^\mathbb{C}$.
\qed
\end{lemma}

Using Lemma~\ref{Zariskilemma} we may assume in the proof of
Theorem~\ref{main} that $U^\mathbb{C}$ coincides with the Zariski
closure of $G$. The proof of the following statement is left to
the reader.

\begin{lemma}
Let $G$ be a connected Zariski-dense Lie subgroup of
$U^\mathbb{C}$ and $N$ a connected complex Lie subgroup of
$U^\mathbb{C}$ which contains $G$. Then every connected normal
complex subgroup of $N$ is a normal subgroup of $U^\mathbb{C}$. In
particular $N$ is normalized by $U^\mathbb{C}$. \qed
\end{lemma}

\begin{remark}
If $U$ is semisimple, then $N=U^\mathbb{C}$.
\end{remark}

An application of the lemma shows that the connected complex Lie
subgroup $G^c$ of $U^\mathbb{C}$ with Lie algebra
$\mathfrak{g}^c:=\mathfrak{g}+i\mathfrak{g}$ is a normal complex
Lie subgroup of $U^\mathbb{C}$ and also that the connected complex
Lie subgroup with Lie algebra $\mathfrak{g}\cap i\mathfrak{g}$ is
normalized by $U^\mathbb{C}$. Consequently $\mathfrak{g}^c$ and
$\mathfrak{g}\cap i\mathfrak{g}$ are ideals in
$\mathfrak{u}^\mathbb{C}$.

On the level of Lie algebras the structure of $G$ is described by the following

\begin{proposition}\label{uZerlegung}
Let $G$ be a connected compatible Zariski-dense Lie subgroup of
$U^\mathbb{C}$. Then there exist ideals $\mathfrak{u}_0$ and
$\mathfrak{u}_1$ in $\mathfrak{u}$, such that
$\mathfrak{g}=\mathfrak{g}_0\oplus\mathfrak{u}_1^\mathbb{C}$ where
$\mathfrak{g}_0\subset\mathfrak{u}_0^\mathbb{C}$ is a compatible
real form of $\mathfrak{u}_0^\mathbb{C}$, i.e.
$\mathfrak{g}_0=\mathfrak{k}_0\oplus\mathfrak{p}_0$ with
$\mathfrak{k}_0\subset\mathfrak{u}_0$, $\mathfrak{p}_0\subset
i\mathfrak{u}_0$ and $\mathfrak{u}_0=\mathfrak{k}_0\oplus
i\mathfrak{p}_0$.
\end{proposition}

\begin{proof}
Let $\mathfrak{u}_1:=\mathfrak{k}\cap i\mathfrak{p}$. Since
$\mathfrak{g}=\mathfrak{k}\oplus\mathfrak{p}\subset\mathfrak{u}
\oplus i\mathfrak{u}$, the equality
$\mathfrak{u}_1^\mathbb{C}=(\mathfrak{k}\cap i\mathfrak{p})\oplus
i(\mathfrak{k}\cap i\mathfrak{p})=\mathfrak{g}\cap i\mathfrak{g}$
holds. Since $G$ is assumed to be connected and Zariski-dense in
$U^\mathbb{C}$, $\mathfrak{u}_1^\mathbb{C}=\mathfrak{g}\cap
i\mathfrak{g}$ is an ideal in $\mathfrak{u}^\mathbb{C}$.
Consequently $\mathfrak{u}_1$ is an ideal in $\mathfrak{u}$.

We fix a $U$-invariant positive inner product $\langle\ ,\
\rangle$ on $\mathfrak{u}$. Note that
$\langle[\xi,\eta_1],\eta_2\rangle=-\langle\eta_1,[\xi,\eta_2]\rangle$
holds for every $\xi, \eta_1, \eta_2\in\mathfrak{u}$.

Let $\mathfrak{k}_0:=\mathfrak{k}\cap\mathfrak{u}_1^\perp$ and
$\mathfrak{p}_0:=\mathfrak{p}\cap i\mathfrak{u}_1^\perp$. Defining
$\mathfrak{u}_0:=\mathfrak{k}_0+i\mathfrak{p}_0$ the sum is direct
by construction and $\mathfrak{u}_0$ is an ideal in $\mathfrak{u}$
since $\langle\ ,\ \rangle$ is assumed to be $U$-invariant,
$\mathfrak{u}_1$ is an ideal in $\mathfrak{u}$ and
$\mathfrak{g}^c$ is an ideal in $\mathfrak{u}^\mathbb{C}$. Setting
$\mathfrak{g}_0:=\mathfrak{k}_0\oplus\mathfrak{p}_0$, it remains
to show that
$\mathfrak{g}=\mathfrak{g}_0\oplus\mathfrak{u}_1^\mathbb{C}$. For
this decompose $\xi\in\mathfrak{g}$ as
$\xi=\xi_1+\xi_{\mathfrak{u}_1^\perp}+\xi_{i\mathfrak{u}_1^\perp}$
with respect to
$\mathfrak{u}^\mathbb{C}=\mathfrak{u}_1^\mathbb{C}\oplus\mathfrak{u}_1^\perp\oplus
i\mathfrak{u}_1^\perp$. Then
$\xi_{\mathfrak{u}_1\perp}+\xi_{i\mathfrak{u}_1^\perp}\in\mathfrak{g}\cap(\mathfrak{u}_1^\perp\oplus
i\mathfrak{u}_1^\perp)=\mathfrak{k}_0\oplus\mathfrak{p}_0=\mathfrak{g}_0$.
\end{proof}

If $U$ is semisimple and $G$ is connected, then
Proposition~\ref{uZerlegung} extends to the group level. Of course
there is also a version for a general compact group $U$ but this
is not needed here.

\begin{proposition}\label{UZerlegung}
Let $U$ be semisimple and let $G$ be a connected compatible
Zariski-dense Lie subgroup of $U^\mathbb{C}$. Then there exist
compact connected Lie subgroups $U_0$ and $U_1$ of $U$ which
centralize each other, such that
\begin{enumerate}
\item $U^\mathbb{C}=U_0^\mathbb{C}\cdot U_1^\mathbb{C}$ and the
intersection $U_0^\mathbb{C}\cap U_1^\mathbb{C}$ is a finite
subgroup of the center of $U$, \item $G=G_0\cdot U_1^\mathbb{C}$,
where $G_0$ is a real form of $U_0^\mathbb{C}$ which is compatible
with the Cartan decomposition, i.e $G_0=K_0\exp(\mathfrak{p}_0)$
where $K_0$ is a Lie subgroup of $U_0$ and $\mathfrak{p}_0$ is a
$K_0$-stable subspace of $i\mathfrak{u}_0$.
\end{enumerate}
\end{proposition}

\begin{proof}
Since $\mathfrak{u}$ is semisimple, it decomposes uniquely into a
sum of simple ideals. Therefore $\mathfrak{u}_0$ and
$\mathfrak{u}_1$ as defined in the proof of
Proposition~\ref{uZerlegung} are semisimple Lie algebras and
consequently the corresponding subgroups $U_0$ and $U_1$ of $U$
are compact. They centralize each other since $\mathfrak{u}_0$ and
$\mathfrak{u}_1$ are ideals. The first statement follows since $G$
is Zariski-dense in $U^\mathbb{C}$ and
$\mathfrak{u}_0\cap\mathfrak{u}_1=\{0\}$.

Defining $K_0:=U_0\cap K$, Proposition~\ref{uZerlegung} implies the
second statement.
\end{proof}

\section{Compatible commutative  subgroups}\label{commutative}

Let $\mathfrak{a}$ be a maximal Lie subalgebra of $\mathfrak{g}$
such that $\mathfrak{a} \subset \mathfrak{p}$. Since
$[\mathfrak{p},\mathfrak{p}]\subset \mathfrak{k}$, the Lie algebra
$\mathfrak{a}$ is commutative. Let $A:=\exp \mathfrak{a}$ denote
the corresponding commutative Lie subgroup of $G$. Note that $A$
is compatible with the Cartan decomposition of $U^\mathbb{C}$ and
let $\mu_{i\mathfrak{a}}\colon Z\to i\mathfrak{a}^*$ be the
restriction of $\mu$ to $i\mathfrak{a}$. The set
$\mathcal{S}_A(\mathcal M_{i\mathfrak{a}})$ where $\mathcal
M_{i\mathfrak{a}}:=\mu_{i\mathfrak{a}}^{-1}(0)$ is open in $Z$
(\cite{HS05}).

For $\xi\in\mathfrak{u}^\mathbb{C}$ let $\xi_Z$ denote the vector
field corresponding to the one parameter group $(t,z)\mapsto \exp
t\xi\cdot z$, i.e.
$\xi_Z(z)=\left.\frac{d}{dt}\right|_0\exp(t\xi)\cdot z$. The Lie
algebra of the isotropy group $A_z=\{a\in A;\ a\cdot z=z\}$ of $A$
at $z\in Z$ is given by $\mathfrak{a}_z:=\{\xi\in \mathfrak{a};\
\xi_Z(z)=0\}$. Let $i\mathfrak{a}_z^\circ=\{\beta\in
i\mathfrak{a}^*;\ \beta\vert i\mathfrak{a}_z=0\}$ denote the
annihilator of $i\mathfrak{a}_z$ in $i\mathfrak{a}^*$.

\begin{proposition} \label{Aconvex}
The image $\mu_{i\mathfrak{a}}(A\cdot z)$ of $A\cdot z$ in
$i\mathfrak{a}^*$ is an open convex subset of the affine subspace
$\mu_{i\mathfrak{a}}(z)+i\mathfrak{a}_z^\circ$ of
$i\mathfrak{a}^*$ and the map $\mu_{i\mathfrak{a}}\colon A\cdot
z\to \mu_{i\mathfrak{a}}(A\cdot z)$ is a diffeomorphism.
\end{proposition}

\begin{proof}
By definition of $\mu_{i\mathfrak{a}}$ we have
\[d\mu_{i\mathfrak{a}}(v)(\xi)=\omega(\xi_Z(y),v)\]
for all $y\in A\cdot z$, $v\in T_y(A\cdot z)$ and $\xi\in
i\mathfrak{a}$. Since $\mathfrak{a}_y=\mathfrak{a}_z$ holds for
all $y\in Az$, we conclude $\mu_{i\mathfrak{a}}(A\cdot
z)\subset\mu_{i\mathfrak{a}}(z)+i\mathfrak{a}_z^\circ$. From the
identification $T_y(A\cdot z)\cong \mathfrak{a}/\mathfrak{a}_y$
and $\ker d(\mu_{i\mathfrak{a}}\vert A\cdot z)=0$ it follows that
$\mu_{i\mathfrak{a}}$ maps $A\cdot z$ locally diffeomorphically
onto an open subset of
$\mu_{i\mathfrak{a}}(z)+i\mathfrak{a}_z^\circ$. Injectivity is
proved in \cite{HS05}, Lemma 5.4.

In order to see that $\mu_{i\mathfrak{a}}(A\cdot z)$ is convex one
may proceed as in \cite{HH96}. In fact if $\exp\,
i\mathfrak{a}\subset U$ would be compact, then the result would
directly follow from the convexity result in \cite{HH96}. For the
convenience of the reader we recall the argument in our situation.

Let $T:=\overline {\exp \, i\mathfrak{a}}$ and denote the Lie
algebra of $T$ by $\mathfrak{t}$. Let $T^\mathbb{C}=\exp\,
(\mathfrak{t}+i\mathfrak{t})$ denote the complexification of the
compact torus $T$. For  $\beta_0$ and $\beta_1$ in
$\mu_{i\mathfrak{a}}(A\cdot z)$ we choose $x_j\in A\cdot z$ such
that $\mu_{i\mathfrak{a}}(x_j)=\beta_j$ holds. Let $\hat
\beta_j:=\mu_\mathfrak{t}(x_j)$. The set $X:=T^\mathbb{C}\cdot z$
is an injectively immersed complex submanifold of $Z$ if we endow
it with the complex structure induced by the natural
identification of $T^\mathbb{C}/T^\mathbb{C}_z$ with
$T^\mathbb{C}\cdot z$. The K\"ahlerian structure on $Z$ restricts
in the sense of pulling back to $X$ and $\mu_\mathfrak{t}|_X$ is a
momentum map on $X$.

Now consider the shifted momentum map
$\hat\mu_j:=\mu_\mathfrak{t}|_X-\hat \beta_j$. From \cite{HS05},
Theorem 10.3, we get the existence of a strictly plurisubharmonic
exhaustion function $\rho_j$ on $X$, such that $\hat\mu_j$ is
associated to $\rho_j$, i.e,
$\hat\mu_j(x)(\xi)=\left.\frac{d}{dt}\right|_0\rho_j(\exp(it\xi)\cdot
x)$ for every $\xi\in\mathfrak{t}$ and $x\in X$. Let $s\in [0,1]$.
We have to show, that
$(1-s)\beta_0+s\beta_1\in\mu_{i\mathfrak{a}}(Az)$. Consider
\[
\hat\mu_s:=(1-s)\hat\mu_0+s\hat\mu_1=\mu_\mathfrak{t}|_X-((1-s)\hat
\beta_0+s\hat \beta_1).
\]
This defines a momentum map on $X$ and it is associated to the
strictly plurisubharmonic exhaustion function
$\rho_s:=(1-s)\rho_0+s\rho_1$.

Since $z$ is contained in the zero-set of the shifted momentum map
$\mu_\mathfrak{t}|_X-\mu_\mathfrak{t}(z)$, the isotropy
$(T^\mathbb{C})_z$ is compatible with the Cartan decomposition
(\cite{HS05}, Lemma 5.5), i.e., equals
$(T_z)^\mathbb{C}=T_z\exp(i\mathfrak{t}_z)$. The orbit $A\cdot
z\cong A/A_z$ is closed in $T^\mathbb{C}/T^\mathbb{C}_z\cong
T^\mathbb{C}\cdot z$. This implies that $\rho_s|_{A\cdot z}$ is an
exhaustion and therefore attains its minimal value at some point
$z_0\in A\cdot z$. Then $\hat\mu_s(z_0)=0$ which implies that
$(1-s)\beta_0+s\beta_1\in\mu_{i\mathfrak{a}}(A\cdot z)$.
\end{proof}

Let $\mathfrak{m}$ be a linear subspace of $\mathfrak{u}$ and
$\eta_{\mathfrak{m}}:=\frac{1}{2}\lVert
\mu_{\mathfrak{m}}\rVert^2$ where $\lVert\
\rVert\colon\mathfrak{u}^*\to\mathbb{R}$ denotes a $U$-invariant
norm function. We identify $\mathfrak{m}$ and $\mathfrak{m}^*$ by
the corresponding inner product $\langle\,,\,\rangle$ on
$\mathfrak{u}$. Since $d\eta_\mathfrak{m}(z)(v)=\langle
d\mu_\mathfrak{m}(z)(v),\mu_\mathfrak{m}(z)\rangle=
\omega((\mu_\mathfrak{m}(z))_Z(z),v)$, we have
$\grad(\eta_\mathfrak{m})(z)=(i\mu_\mathfrak{m}(z))_Z(z)$ with
respect to the Riemannian structure induced by $\omega$ on smooth
$U^\mathbb{C}$-stable submanifolds of $Z$. Note that $z\mapsto
(i\mu_{\mathfrak{m}}(z))_Z(z)$ is a globally defined vector field
on $Z$. By abuse of notation, we denote it by
$\grad(\eta_\mathfrak{m})$ and call it the gradient vector field
of $\eta_{\mathfrak{m}}$. The set
$\mathcal{C}_{\mathfrak{m}}:=\{z\in Z;\
\grad(\eta_\mathfrak{m})(z)=0\}$ is called the set of critical
points of $\eta_\mathfrak{m}$.

For $\mathfrak{m}=i\mathfrak{a}$ we have

\begin{lemma}\label{acritical}
\( \mathcal{S}_A(\mathcal M_{i\mathfrak{a}})\cap \mathcal
C_{i\mathfrak{a}}=\mathcal M_{i\mathfrak{a}}\,. \)
\end{lemma}

\begin{proof}
The inclusion $\mathcal
M_{i\mathfrak{a}}\subset\mathcal{S}_A(\mathcal
M_{i\mathfrak{a}})\cap\mathcal C_{i\mathfrak{a}}$ is trivial.

Let $z\in\mathcal{S}_A(\mathcal M_{i\mathfrak{a}})\cap \mathcal
C_{i\mathfrak{a}}$. From Proposition~\ref{Aconvex} we know that
$\mu_{i\mathfrak{a}}(Az)\subset
\mu_{i\mathfrak{a}}(z)+i\mathfrak{a}_z^\circ$. Since
$z\in\mathcal{S}_A(\mathcal M_{i\mathfrak{a}})$, we have
$0\in\mu_{i\mathfrak{a}}(\overline{Az})\subset
\mu_{i\mathfrak{a}}(z)+i\mathfrak{a}_z^\circ$, so
$\mu_{i\mathfrak{a}}(Az)\subset i\mathfrak{a}_z^\circ$. Then
$0=\grad(\eta_\mathfrak{a})(z)=(i\mu_{i\mathfrak{a}}(z))_Z(z)$,
i.e. $i\mu_{i\mathfrak{a}}(z)\in \mathfrak{a}_z^*$ implies
$\mu_{i\mathfrak{a}}(z)=0$.
\end{proof}

\noindent For $r>0$ let
$\Delta_r(\eta_{i\mathfrak{a}}):=\eta_{i\mathfrak{a}}^{-1}([0,
r))$.

\begin{proposition} \label{Acompactneighborhood}
Let $C$ be a compact subset of $\mathcal M_{i\mathfrak{a}}$ and
let $W$ be an open neighborhood of $C$ in $Z$. Then there exist an
$A$-stable open set $\Omega$ in $Z$ and an $r>0$ such that
$C\subset \Omega\cap\Delta_r(\eta_{i\mathfrak{a}})\subset W$
holds.
\end{proposition}

\begin{proof}
In the proof we use the results in \cite{HS05} freely.

Since $A$ is commutative $\mathcal{S}_A(\mathcal
M_{i\mathfrak{a}})$ is open in $Z$ and we  may assume that
$Z=\mathcal{S}_A(\mathcal M_{i\mathfrak{a}})$ holds. In particular
the topological Hilbert quotient $\pi\colon Z\to Z/\!\!/ A$
exists. Since the maximal compact subgroup of $A$ is trivial the
map $\pi\circ\imath\colon \mathcal M_{i\mathfrak{a}}\to Z/\!\!/ A$
where $\imath \colon\mathcal M_{i\mathfrak{a}} \to Z$ is the
inclusion is a homeomorphism. If we identify $Z/\!\!/ A$ with
$\mathcal M_{i\mathfrak{a}}$ the quotient map $\pi\colon
Z\to\mathcal M_{i\mathfrak{a}}$ is given by $z\mapsto
\overline{A\cdot z}\cap\mathcal M_{i\mathfrak{a}}$ and the
$\pi$-fiber trough $q\in\mathcal M_{i\mathfrak{a}}$ is given by
$\pi^{-1}(q)=\mathcal{S}_A(q)=\{z\in Z;\ q\in \overline{A\cdot
z\cdot}\}$.

We may assume that $W$ is relatively compact. Let $V_{\mathcal
M_{i\mathfrak{a}}}$ be a relatively compact open neighborhood of
$C$ in $W\cap\mathcal M_{i\mathfrak{a}}$ and let $r^*$ be the
minimal value of $\eta_{i\mathfrak{a}}$ restricted to the compact
set $\partial W\cap\mathcal{S}_A(\overline{V_{\mathcal
M_{i\mathfrak{a}}}})$. Note that $r^*>0$. For $r<r^*$ let
$V_r:=W\cap\Delta_r(\eta_{i\mathfrak{a}})\cap\mathcal{S}_A(V_{\mathcal
M_{i\mathfrak{a}}})$.
\\
We claim that $V_r=A\cdot V_r\cap\Delta_r(\eta_{i\mathfrak{a}})$
holds. For this we have to show that $A\cdot z
\cap\Delta_r(\eta_{i\mathfrak{a}})\subset V_r$ for every $z\in
V_r$.

Proposition~\ref{Aconvex} implies that $A\cdot z
\cap\Delta_r(\eta_{i\mathfrak{a}})$ is connected. By definition of
$r$, it does not intersect the boundary of $W$. Since $z\in W$, we
conclude $A\cdot z \cap\Delta_r(\eta_{i\mathfrak{a}})\subset W$
and therefore $A\cdot z \cap\Delta_r(\eta_{i\mathfrak{a}})\subset
V_r$.

Setting $\Omega:=A\cdot V_r$, the proposition follows.
\end{proof}

\section{Neighborhoods of the zero fiber}

Let $G=K\exp\mathfrak{p}$ be a closed compatible Lie subgroup of
$U^\mathbb{C}$ and $\mathfrak{a}$ a maximal subalgebra of
$\mathfrak{p}$. We have $G=KAK$ or equivalently
$\mathfrak{p}=K\cdot\mathfrak{a}$ with respect to the adjoint
action of $K$ on $\mathfrak{p}$.

Recall that $\mathcal C_{i\mathfrak{p}}$ is by definition the set
of zeros of the vector field
$z\mapsto(i\mu_{i\mathfrak{p}}(z))_Z(z)$ where $i\mathfrak{p}$ and
$i\mathfrak{p}^*$ are identified by a $U$-invariant positive inner
product on $\mathfrak{u}$. A result analogous to
Lemma~\ref{acritical} holds for $G$. More precisely, for
$X:=\bigcap_{k\in K}\, k\cdot\mathcal{S}_A(\mathcal
M_{i\mathfrak{a}})$ we have

\begin{lemma}\label{gcritical}
\( X\cap \mathcal{C}_{i\mathfrak{p}}=\mathcal M_{i\mathfrak{p}}\,.
\)
\end{lemma}

\begin{proof}
We have $\mathcal M_{i\mathfrak{p}}\subset X\cap
\mathcal{C}_{i\mathfrak{p}}$ since $\mathcal M_{i\mathfrak{p}}$ is
$K$-invariant.

Let $z\in X\cap \mathcal{C}_{i\mathfrak{p}}$. There exists a $k\in
K$ with $\mu_{i\mathfrak{p}}(k\cdot
z)=k\cdot\mu_{i\mathfrak{p}}(z)\in i\mathfrak{a}^*$. Then
$z':=k\cdot z\in X\subset\mathcal{S}_A(\mathcal
M_{i\mathfrak{a}})$ since $X$ is $K$-invariant and $z'\in\mathcal
C_{i\mathfrak{a}}$ since
$\mu_{i\mathfrak{p}}(z')=\mu_{i\mathfrak{a}}(z')$ and $z'\in
\mathcal{C}_{i\mathfrak{p}}$. Lemma~\ref{acritical} implies
$z'\in\mathcal M_{i\mathfrak{a}}$. We conclude $z'\in
X\cap\mathcal M_{i\mathfrak{p}}$. Invariance of $\mathcal
M_{i\mathfrak{p}}$ with respect to $K$ implies $z\in X\cap\mathcal
M_{i\mathfrak{p}}$.
\end{proof}

For $r>0$ let
$\Delta_r(\eta_{i\mathfrak{p}}):=\eta_{i\mathfrak{p}}^{-1}([0,
r))$. Now we give a generalization of
Proposition~\ref{Acompactneighborhood} to the non-commutative
case.

\begin{proposition}\label{Gcompactneighborhood}
Let $C$ be a compact $K$-stable subset of $\mathcal
M_{i\mathfrak{p}}$ and let $W$ be an open neighborhood of $C$ in
$Z$. Then there exist a $G$-stable open set $\Omega$ in $Z$ and an
$r>0$ such that $C\subset
\Omega\cap\Delta_r(\eta_{i\mathfrak{p}})\subset W$ holds.
\end{proposition}

\begin{proof}
We may assume that $W$ is $K$-stable. From
Proposition~\ref{Acompactneighborhood} we get the existence of an
$A$-invariant open neighborhood $\Omega_A$ of $C$ and an $r>0$,
such that $\Omega_A\cap\Delta_r(\eta_{i\mathfrak{a}})\subset W$
holds. Let $V$ be a $K$-stable stable neighborhood of $C$ in
$\Omega_A$. Then $A\cdot
V\cap\Delta_r(\eta_{i\mathfrak{a}})\subset W$. Let $\Omega:=G\cdot
V$. Since $\eta_{i\mathfrak{a}}(z)\leq\eta_{i\mathfrak{p}}(z)$ we
conclude
\begin{align*}
&\Omega\cap\Delta_r(\eta_{i\mathfrak{p}})=K\cdot A\cdot
V\cap\Delta_r(\eta_{i \mathfrak{p}})
=K\cdot(A\cdot V\cap\Delta_r(\eta_{i\mathfrak{p}}))\\
 &\subset K\cdot(A\cdot
V\cap\Delta_r(\eta_{i\mathfrak{a}}))\subset K\cdot W= W\,.
\end{align*}
\end{proof}

\begin{remark}
The proof uses
$\Delta_r(\eta_{i\mathfrak{p}})\subset\Delta_r(\eta_{i\mathfrak{a}})$.
More precisely it can be shown that
$\Delta_r(\eta_{i\mathfrak{p}})=\bigcap_{k\in K}\,
k\cdot\Delta_r(\eta_{i\mathfrak{a}})$ holds.
\end{remark}

\section{The set of semistable points is open}

Recall that $\grad(\eta_{i\mathfrak{p}})$ denotes the vector field
$z\mapsto (i\mu_{i\mathfrak{p}}(z))_Z(z)$. We have
$\grad(\eta_{i\mathfrak{p}})(z)\in T_z(G\cdot z)$, i.e., it is
tangent to the $G$-orbits. Since the momentum map is assumed to be
smooth on submanifolds of $Z$, the gradient vector field is smooth
on $U^\mathbb{C}$-stable submanifolds. Let $\psi_z\colon I_z\to
Z$, $t\mapsto \psi_z(t)$ denote integral curve of
$\grad(\eta_{i\mathfrak{p}})$ with $\psi_z(0)=z$. Here we assume
the interval $I_z$ to be maximal.

\begin{theorem}
The set $\mathcal{S}_G(\mathcal M_{i\mathfrak{p}})$ is open in
$Z$.
\end{theorem}

\begin{proof}
We may assume that $G$ is closed in $U^\mathbb{C}$. Let
$x\in\mathcal M_{i\mathfrak{p}}$. It is sufficient to show that an
open neighborhood $V$ of $x$ is contained in
$\mathcal{S}_G(\mathcal M_{i\mathfrak{p}})$. By
Proposition~\ref{Gcompactneighborhood} there exists a $G$-stable
open neighborhood $\Omega$ of $x$, such that for some $r>0$ the
set $\overline{\Omega\cap\Delta_r(\eta_{i\mathfrak{p}})}$ is
compact and contained in $X=\bigcap_{k\in K}\,
k\cdot\mathcal{S}_A(\mathcal M_{i\mathfrak{a}})$. Let
$V:=\Omega\cap\Delta_r(\eta_{i\mathfrak{p}})$. By definition of
$\psi$ we have $\psi_z(t)\in G\cdot
z\cap\overline{\Delta_r(\eta_{i\mathfrak{p}})}$ for all
$z\in\overline{\Delta_r(\eta_{i\mathfrak{p}})}$ and all $t\le 0$
in the domain of definition of $\psi_z$. Let $s_z\in [-\infty,0)$
be minimal with $(s_z,0]\subset I_z$. Note that $s_z=-\infty$ for
$z\in \overline V$ if $Z$ is smooth.

Let $Z=R_0\cup R_1\cup\ldots\cup R_k$ be the stratification of $Z$
into smooth submanifolds, i.e. $R_0$ is the set of smooth  points
in $Z$ and $R_{l+1}$ is the set of smooth  points in
$Z\setminus(R_0\cup \dotsb \cup R_l)$. Let $z_0\in V$ and let
$t_n\leq 0$ be a decreasing sequence which converges to $s_{z_0}$.
After possibly replacing $(t_n)$ by a subsequence, the limit
$z_1:=\lim_{n\to\infty}\psi_{z_0}(t_n)$ exists. By induction we
get a sequence $(z_n)$. If $z_n\in R_j$ and $s_{z_n}\neq -\infty$,
then $z_{n+1}\in R_{j+l}$ for an $l>0$ since $R_j$ is smooth and
$U^\mathbb{C}$-invariant. Therefore there exists an $n_0$ with
$s_{z_{n_0}}=-\infty$. To simplify notation, we assume $n_0=0$.

Let $t_0<0$. Then
$\psi_{z_1}(t_0)=\lim_{n\to\infty}\psi_{z_0}(t_n+t_0)$ and
consequently
\[\eta_{i\mathfrak{p}}(\psi_{z_1}(t_0))=\lim_{n\to\infty}\eta_{i\mathfrak{p}}(\psi_{z_0}(t_n+t_0))=\lim_{n\to\infty}\eta_{i\mathfrak{p}}(\psi_{z_0}(t_n))=\eta_{i\mathfrak{p}}(z_1).\]
Hence $\psi_{z_1}(t_0)=z_1$ and $z_1\in\mathcal C_{i\mathfrak{p}}$.\\
Lemma~\ref{gcritical} implies that  $z_1\in\mathcal
C_{i\mathfrak{p}}\cap\overline V\cap\overline{G\cdot
z_0}\subset\mathcal C_{i\mathfrak{p}}\cap X\cap\overline{G\cdot
z_0}=\mathcal M_{i\mathfrak{p}}\cap\overline{G\cdot z_0}$. This
shows $V\subset\mathcal{S}_G(\mathcal M_{i\mathfrak{p}})$.

\end{proof}

\section{Shifting}\label{shifting}

In this section we consider the case of a general
$\mu_{i\mathfrak{p}}$-fiber $\mathcal M_{i\mathfrak{p}}(\beta)$.
The first step is to shift the momentum map relative to the
coadjoint orbit $U\cdot \beta$. The relevant properties are
summarized in

\begin{lemma}\label{coadjoint}
Let $\beta\in \mathfrak{u}^*$ and let $O:=U\cdot\beta$ be the
coadjoint orbit of $\beta$. Moreover let $\xi\in \mathfrak{u}$
denote the dual vector of $\beta$ with respect to a $U$-invariant
negative definite inner product on $\mathfrak{u}$. Let $Q:= \{g\in
U^\mathbb{C};\ \lim_{t\to-\infty}\exp(it\xi)g\exp(-it\xi)\text{
exists in } U^\mathbb{C}\}$. Then the following holds
\begin{enumerate}
\item $Q$ is a parabolic subgroup of $U^\mathbb{C}$ with Lie
algebra $\mathfrak{q}=\mathfrak{z}^\mathbb{C}(\xi)\oplus
\mathfrak{r}$ where $\mathfrak{z}^\mathbb{C}(\xi)$ denotes the
centralizer of $\xi$ in $\mathfrak{u}^\mathbb{C}$ and
$\mathfrak{r}$ denotes the sum of eigenspaces of $\ad(-i\xi)$ with
negative eigenvalues.
\item The $U$-equivariant map $\imath\colon
O\to U^\mathbb{C}/Q$, $u\cdot\beta\mapsto uQ$
      is a real analytic isomorphism.
\item  The $U$-invariant symplectic structure
$\omega_O(\zeta_O(\alpha),\eta_O(\alpha)):=-\alpha([\zeta,\eta])$
for $\alpha\in O$ and $\zeta,\eta\in\mathfrak{u}$ is a K\"ahlerian
structure with respect to the complex structure on $O$ induced by
$\imath$. The map $O\to \mathfrak{u}^*$, $\eta\mapsto -\eta$
defines a momentum map on $O$.\qed
\end{enumerate}
\end{lemma}

Lemma~\ref{coadjoint} says that $O=U\cdot \beta$ is a complex
manifold such that $U^\mathbb{C}$ acts holomorphically and
transitively, i.e., $U\cdot \beta=U^\mathbb{C}\cdot \beta$ holds.
The following observation is crucial for our application of the
shifting procedure.

\begin{proposition}
$G\cdot\beta=K\cdot \beta$ for $\beta\in i\mathfrak{p}^*$.
\end{proposition}

\begin{proof}
As a first step we replace $U^\mathbb{C}$ by its image
$\Ad(U)^\mathbb{C}$ under the adjoint representation. The image
$\Ad(G)$ is compatible with the Cartan decomposition of
$\Ad(U)^\mathbb{C}$. Then $U$ is semisimple.

Assuming in addition that $G$ is Zariski-dense in $U^\mathbb{C}$
(Lemma~\ref{Zariskilemma}) and connected, we may apply
Proposition~\ref{UZerlegung}. It yields that $U$ is the product of
two compact subgroups $U_0$ and $U_1$ which centralize each other,
such that $G=G_0\cdot U_1^\mathbb{C}$. We decompose $\beta\in
i\mathfrak{p}^*$ with respect to this decomposition, i.e.
$\beta=\beta_0+\beta_1$ where $\beta_0\in
i\mathfrak{p}_0^*=i(\mathfrak{p}\cap i\mathfrak{u}_0)^*$ and
$\beta_1\in i\mathfrak{u}_1^*$. The coadjoint orbit $O=U\cdot
\beta$ can be biholomorphically and $U$-equivariantly identified
with $U_0\cdot\beta_0\times U_1\cdot\beta_1$. We have
$U_1^\mathbb{C}\cdot\beta_1=U_1\cdot\beta_1$ by construction, so
we may assume $U=U_0$.

Let $\sigma\colon \mathfrak{u}^\mathbb{C}\to
\mathfrak{u}^\mathbb{C}$ be the anti-holomorphic involution of Lie
algebras defined by $\sigma\vert\mathfrak{g}=\id_{\mathfrak{g}}$
and $\sigma\vert i\mathfrak{g}=-\id_{i\mathfrak{g}}$. Since
$\beta\in i\mathfrak{p}^*$, the Lie algebra $\mathfrak{q}$ of $Q$
is stable with respect to $\sigma$. Therefore it defines an
antiholomorphic involution on the tangent space $T_{\beta}
(U\cdot\beta)\cong\mathfrak{u}^\mathbb{C}/\mathfrak{q}$. For a
subspace $\mathfrak{m}\subset \mathfrak{u}^\mathbb{C}$ define
$\mathfrak{m}\cdot\beta:=\{\xi_O(\beta);\xi\in\mathfrak{m}\}$. We
have $T_\beta
(U\cdot\beta)=\mathfrak{u}^\mathbb{C}\cdot\beta=\mathfrak{k}\cdot
\beta+i\mathfrak{p}\cdot\beta=i\mathfrak{k}\cdot\beta+\mathfrak{p}\cdot\beta$.
Since
$\sigma\vert(\mathfrak{g}\cdot\beta)=\id_{\mathfrak{g}\cdot\beta}$
and $\sigma\vert
(i\mathfrak{g}\cdot\beta)=-\id_{i\mathfrak{g}\cdot\beta}$, it
follows that $\mathfrak{k}\cdot\beta=\mathfrak{p}\cdot\beta$.
Therefore $T_\beta (G\cdot\beta)=T_\beta (K\cdot\beta)$ and we
conclude $G\cdot\beta= K\cdot\beta$.
\end{proof}

Now we can prove the main result of this paper.
For the convenience of the reader we repeat the
statement.

\begin{theorem11}
For every $\beta\in i\mathfrak{p}^*$ the set
$\mathcal{S}_G(\mathcal M_{i\mathfrak{p}}(\beta))$ is open in $Z$.
\end{theorem11}

\begin{proof}
The K\"ahlerian structures on $Z$ and $O$ induce a K\"ahlerian
structure on the product $Z\times O$. It follows that
$\hat\mu\colon Z\times O\to\mathfrak{u}^*$,
$\hat\mu(z,\alpha)=\mu(z)-\alpha$ is a $U$-equivariant momentum
map on $Z\times O$ with respect to this product K\"ahlerian
structure and the diagonal action of $U$. Moreover
$\mathcal{S}_G(\widehat{\mathcal M}_{i\mathfrak{p}})$, where
$\widehat{\mathcal
M}_{i\mathfrak{p}}:=\hat\mu_{i\mathfrak{p}}^{-1}(0)$  is open in
$Z\times O$ and therefore its intersection $Z(\beta)\subset
Z\times O$  with $Z\times G\cdot \beta$ is open as well. The
projection $p\colon Z\times G\cdot \beta\to Z$ is an open map.
Since $G\cdot\beta$ is a $K$-orbit, $p$ maps the open $G$-stable
set $Z(\beta)$ onto $\mathcal{S}_G(\mathcal
M_{i\mathfrak{p}}(\beta))$.
\end{proof}


\newcommand{\noopsort}[1]{} \newcommand{\printfirst}[2]{#1}
\newcommand{\singleletter}[1]{#1} \newcommand{\switchargs}[2]{#2#1}
\providecommand{\bysame}{\leavevmode\hbox to3em{\hrulefill}\thinspace}

\end{document}